\documentclass[11pt,reqno]{amsart}
\usepackage{amsmath,amssymb,amsbsy,amsfonts,amsthm,latexsym,
            amsopn,amstext,amsxtra,euscript,amscd,amsthm}

\newtheorem{thm}{Theorem}
\newtheorem{theorem}[thm]{Theorem}

\newtheorem{conj}{Conjecture}

\theoremstyle{definition}

\newtheorem*{rem*}{Remark}

\def\\{\cr}
\def\({\left(}
\def\){\right)}
\def\[{\left[}
\def\]{\right]}
\def\<{\langle}
\def\>{\rangle}

\def\ge{\geqslant}
\def\le{\leqslant}

\begin{document}

\title[Realizable subsequences of linearly recurrent sequences]{On (almost) realizable subsequences of linearly recurrent sequences}

\author{Florian Luca}
\address{School of Mathematics\\
University of the Witwatersrand\\
Private Bag 3, Wits 2050, South Africa}
\address{
Research Group in Algebric Structures and Applications \\
King Abdulaziz University, Jeddah 22254, Saudi Arabia}
\address{ Max Planck Institute for Software Systems\\
Saarland Information Campus E1 5, 66123, Saarbr\"ucken Germany}
\address{Centro de Ciencias Matem\'{a}ticas UNAM \\
Morelia, Mexico}
\email{Florian.Luca@wits.ac.za}

\author{Tom Ward}
\email{tbward@gmail.com}

\date{\today}


\pagenumbering{arabic}

\maketitle

\begin{abstract}
In this note we show that if~$(u_n)_{n\ge 1}$ is a simple linearly recurrent sequence of integers whose minimal recurrence of order~$k$ involves
only positive coefficients that has positive initial terms,
then~$(Mu_{n^s})_{n\ge 1}$ is the sequence of periodic point
counts for some map for a suitable positive
integer~$M$ and~$s$ any sufficiently large multiple
of~$k!$. This extends a result of Moss and Ward [{\it The Fibonacci Quarterly} {\bf 60} (2022), 40--47] who proved the result for
the Fibonacci sequence.
\end{abstract}

\section{Introduction}

A sequence of nonnegative integers~$(a_n)_{n\ge 1}$
is called \emph{realizable} if there is some set~$X$ and a
map~$T\colon X\to X$ such that
\begin{equation*}
a_n={\text{\rm Fix}}(T^n)=\#\{x\in X\mid T^n(x)=x\},
\end{equation*}
and is called \emph{almost realizable} if
it is realizable after multiplication by a constant.

A simple example of realizability is the shift map~$T\colon(x_n)_{n\in {\mathbb Z}}\mapsto (x_{n+1})_{n\in {\mathbb Z}}$
on the {\it golden mean} shift space
$$
X=\{{\bf x}=(x_n)_{n\in {\mathbb Z}}\in \{0,1\}^{\mathbb Z}
\mid(x_k,x_{k+1})\ne (1,1)~\text{for~all}~k\in\mathbb Z\}.
$$
This has
\begin{equation}\label{equationLucasExample}
{\text{\rm Fix}}(T^n)=\#\{{\bf x}\in X\mid T^n({\bf x})={\bf x}\}={\text{\rm Trace}}\left(\begin{matrix} 1 & 1\\ 1 & 0\end{matrix}\right)^n=L_n
\end{equation}
where $(L_n)_{n\ge1}=(1,3,4,7,\ldots)$ is the Lucas companion of the Fibonacci sequence.

A simple example of almost realizability is given by the
sequence
$$
(2^{n-1})_{n\ge1}=(1,2,4,\dots)
$$
of powers of~$2$. Clearly this is not realizable, since
a map~$T$ witnessing this would have~${\rm{Fix}}(T)=1$,
and hence must have~${\rm{Fix}}(T^2)$ odd.
However the shift map as above on the
\emph{full~$2$-shift}~$X=\{0,1\}^{\mathbb{Z}}$
has~${\rm{Fix}}(T^n)=2^n$ for all~$n\ge1$,
showing that the sequence becomes realizable after
multiplication by~$2$.

We refer to~\cite{MT}
for this and other examples,
and the references therein for background on this concept.
Realizable sequences can be characterized
in algebraic terms as follows. The sequence~$(a_n)_{n\ge 1}$
of non-negative integers is realizable
if and only if it satisfies the following two conditions:
\begin{itemize}
\item[(D)] $\displaystyle\sum_{d\vert n}\mu(n/d)a_d\equiv 0\pmod n$ for all $n\in{\mathbb N}$, and
\item[(S)] $\displaystyle\sum_{d\vert n} \mu(n/d)a_d\ge 0$ for all $n\in {\mathbb N}$.
\end{itemize}
Here~$\mu$ denotes the classical M\"obius function.
We call~(D) the \emph{Dold condition} and~(S) the \emph{sign condition}.
The equivalence is clear, because the sum arising in~(D)
and in~(S) is the number of points that lie
on a closed orbit of minimal length~$n$ under iteration of
a map that witnesses realizability of~$(a_n)_{n\ge1}$.

Moss~\cite{pm} showed that~$(5F_{n^2})_{n\ge1}$ is realizable,
and Moss and Ward~\cite{MT} extended this
to show that~$(5F_{n^{2k}})_{n\ge 1}$ is realizable
for~$k\ge1$
while~$(MF_{n^{2k+1}})$ is not realizable for any
choice of~$M=M_k\ge1$ for any~$k\ge0$, where
we write~$(F_n)_{n\ge 1}$
for the Fibonacci sequence~$(1,1,2,\dots)$. These arguments use
congruence properties specific to the Fibonacci sequence,
which makes one wonder to what extent such a result
can be generalised to other linearly recurrent sequences.

Here we generalize the above result to simple linearly recurrent sequences
satisfying some positivity conditions. Our result is quite general, and its proof uses elementary algebraic number theory rather than congruences
specific to a given sequence.

Let~$(u_n)_{n\ge 1}$ be a linearly recurrent sequence of integers of order~$k$. That is, it satisfies a recurrence relation of the form
$$
u_{n+k}=a_1u_{n+k-1}+\cdots+a_k u_{n}
$$
for all~$n\ge 1$,
where~$a_1,\ldots,a_k,u_1,\ldots,u_k$ are all integers
(this is a harmless assumption for
integer linear recurrences by
Fatou's lemma~\cite[p.~369]{MR1555035}). We assume that the recurrence is minimal, so in particular~$a_k\ne 0$.
We ask if there is a polynomial~$f\in {\mathbb Z}[X]$ and a positive integer~$M$
with the property that~$(Mu_{f(n)})_{n\ge 0}$ satisfies~(D) and~(S).
Let
\begin{equation}\label{eq:F}
F(X)=X^k-a_1X^{k-1}-\cdots-a_k,
\end{equation}
the \emph{characteristic polynomial} of
the sequence~$(u_n)_{n\ge 1}$. For background and relevant properties of linearly recurrent sequences the reader is invited to consult the
monograph~\cite{MR1990179}.
Let~${\mathbb K}$ be the splitting field of~$F$
and~${\mathcal O}_{\mathbb K}$ be
its ring of integers. Let $\Delta({\mathbb K})$ be the discriminant of ${\mathbb K}$ and $\Delta(F)$ be the discriminant of ${\mathbb F}$. Let~$G$ be the Galois
group of~${\mathbb K}$ over~${\mathbb Q}$, and
let~$e:=e(G)$ be the exponent of~$G$ and $N$ be the order of $G$.

\begin{theorem}\label{thm:1}
Assume that~$F$ has only simple zeros.
\begin{itemize}
\item[{\rm{(i)}}] The sequence~$(Mu_{n^s})_{n\ge 1}$ satisfies~{\rm{(D)}}
if
\begin{equation}\label{eq:Mu}
\left.\begin{aligned}
M&={\text{\rm lcm}}[\Delta({\mathbb K}), \Delta(F)]\mbox{ and}\\
s&\equiv 0\pmod {e(G)}\mbox{ with }s\ge N.
\end{aligned}
\right\}
\end{equation}
\item[{\rm{(ii)}}]  Assume in addition that~$a_i\ge 0$ for $i=1,\ldots,k$
and~$a_k\ne 0$, that
$$
(a_1,\ldots,a_k)\ne (0,0,\ldots,1),
$$
and that~$u_i\ge 1$ for all~$i\in \{1,2,\ldots,k\}$.
Then the sequence~$(Mu_{n^s})_{n\ge 1}$ satisfies~{\rm(S)} whenever~$\ell\ge \ell_0$ is a sufficiently large
number which can be computed in terms of the sequence~$(u_n)_{n\ge 1}$.
\end{itemize}
\end{theorem}

The somewhat strange condition~(i) can be explained as follows.
The condition~$\Delta(F)\vert M$ is needed to ensure that in the Binet formula
for the general term of~$Mu_n$, the summands involved are algebraic integers.
On the other hand, the additional conditions that~$\Delta({\mathbb K})\vert M$ together with the conditions on~$s$
are sufficient to ensure that the Dold condition~(D) is satisfied.
In particular, if~$F$ is irreducible, then~$\Delta({\mathbb K})\vert\Delta(F)$,
so it suffices that~$M$ is a multiple of~$\Delta(F)$ but this is not
true for reducible polynomials as the example~$F(X)=(X^2-2)(X^3-5)$ shows.
In this case~${\mathbb K}={\mathbb Q}({\sqrt{2}}, {\sqrt[3]{5}},{\sqrt{-3}})$ has~$\Delta(F)=-2^3\cdot 3^3\cdot 5^2\cdot 17^2$,
which is not a multiple of~$\Delta({\mathbb K})=2^{18} 3^{14} 5^8$.

\begin{rem*}
(a) We exclude
a periodic sequence of period~$k\ge2$ and minimal polynomial~$X^k-1$,
since for it the theorem is not true. Indeed,
consider the simplest case when~$k$ is prime. Then condition~(S)
requires that~$u_{k^{s}}-u_1\ge 0$.
On the other hand, since~$k^s\equiv k\pmod k$,
we must have~$u_k\ge u_1$.
It follows that~$u_1,\ldots,u_k$ cannot be chosen to be
arbitrary positive integers. For composite~$k$,
there must be inequalities satisfied between
the
values of~$u_i$ for~$i\in \{1,\ldots,k\}$ whose residue classes
modulo~$k$ are idempotents in the ring~${\mathbb Z}/k{\mathbb Z}$.

\noindent(b) We
require the minimal polynomial to only have simple roots,
for otherwise~(D) may be false.
For example, the sequence
defined by~$u_n=n$ for all~$n\ge 1$,
satisfies the linear recurrence~$u_{n+2}=2u_{n+1}-u_n$
with minimal polynomial~$(X-1)^2$.
Then condition~(D) for~$(Mu_{n^s})$
implies that for a prime~$p$ we must
have~$p\vert M(p^s-1)$,
and for given positive integers~$M$ and~$s$
this can only hold for the finitely many primes~$p$
dividing~$M$. Hence, for the above sequence, the Dold quotients are
rational numbers whose denominators are
divisible by arbitrarily large primes.
Indeed, similar arguments may be used to show
that if~$(f(n))_{n\ge1}$ is realizable
with~$f\in\mathbb{Z}[n]$,
then~$f$ is a constant~\cite[Lem.~2.4]{MR1873399}.

\noindent(c) A different (arguably more natural)
question is to ask when a linear recurrence sequence
itself satisfies~(D) without
multiplication by a factor or passing to a subsequence.
Minton~\cite{MR3195758} showed
that---up to a finite multiplying factor---this
is possible if and only if the sequence is a linear
combination of traces of powers of algebraic numbers.
From this perspective~\eqref{equationLucasExample} is a
manifestation of the fact that the only linearly
recurrent sequences satisfying the Fibonacci
recurrence~$u_{n+2}=u_{n+1}+u_n$ for~$n\ge1$ which have
this property must have~$u_2=3u_1$ (and hence
must be a multiple of the Lucas sequence),
a special case shown earlier in~\cite{MR1866354}
again using congruences specific to the Fibonacci
sequence.
\end{rem*}

Returning to the Fibonacci sequence
where this phenomena was first observed,
the sequence~$(F_{n})_{n\ge 1}$
has~$k=2$,~$a_1=a_2=1$,~$F_1=F_2=1>0$,
and~$F(X)=X^2-X-1$. Further, ${\mathbb K}={\mathbb Q}[{\sqrt{5}}]$.
Thus,~$\Delta(F)=\Delta({\mathbb K})=5$
and~$G={\mathbb Z}/2{\mathbb Z}$ so
$e(G)=N=2$. Thus, $s=2$ satisfies that $s$ is a multiple of $e(G)$ and $s\ge N$.
We shall justify that in this case we can take~$\ell_0=1$, recovering
the result of~\cite{MT} precisely. In fact,
we prove it in a more general setting. Let~$(F_n)^{(k)})_{n\ge -(k-2)}$ be
the~$k$-generalized Fibonacci sequence
satisfying the recurrence
relation~$F_{n+k}^{(k)}=F_{n+k-1}^{(k)}+\cdots+F_n^{(k)}$
for~$n\ge 2-k$
with initial values~$F_i^{(k)}=0$ for~$i=2-k,~3-k,\ldots,-1,~0$ and~$F_1^{(k)}=1$. Wolfram~\cite{MR1622060} conjectured that~$G=S_k$ is
the full symmetric group on $k$ letters,
and this is known to be so when~$k$ is even,
when~$k$ is small, or
when~$k$ is prime (see, for example,
the work of Martin~\cite{MR2043329}).
We therefore take~$N_k:=k!$ and this is a multiple of $e(G)$ and at least as large as $N$.

\begin{theorem}
\label{thm:2}
For~$k\ge 2$
we can take~$s=N_k\ell$ for any $\ell\ge 1$
for the sequence~$(F_{n}^{(k)})_{n\ge 1-(k-2)}$.
\end{theorem}

At the end of their paper~\cite{MT},
Moss and Ward propose the following conjecture.

\begin{conj}
\label{Conj}
Let~$P,~Q\in {\mathbb Z}$
and~$u_{n+2}=Pu_{n+1}+Qu_n$ for~$n\ge1$,
with initial conditions~$u_0=0,~u_1=1$. Then~$((P^2-4Q)u_{n^2})_{n\ge 1}$
satisfies~{\rm{(D)}}.
\end{conj}

This almost follows from Theorem~\ref{thm:1},
except that there are some additional hypotheses
like the fact that the recurrence must be of minimal
order~$k=2$ and that~$(P,Q)\ne (0,1)$, in order to apply the theorem.
We therefore supply a proof of the following result.

\begin{theorem}
\label{thm:3}
\label{conj}
Conjecture~\ref{Conj} holds.
\end{theorem}

\section{The proof of Theorem~\ref{thm:1}}

To start with, let
$$
F(X)=\prod_{i=1}^k (X-\lambda_i),
$$
so that we have the (generalized) Binet formula
$$
u_n=\sum_{i=1}^k c_i\lambda_i^n
$$
for all~$n\ge1$
for coefficients~$c_1,\ldots,c_k$
determined from~$u_1,\ldots,u_k$
by solving a linear system of~$k$ equations in~$k$ unknowns
whose matrix is Vandermonde on~$\lambda_1,\ldots,\lambda_k$.
We write~${\mathbb K}:={\mathbb Q}(\lambda_1,\ldots,\lambda_k)$.
Then~$c_1,\ldots,c_k$ are algebraic numbers in~${\mathbb K}$
having
the Vandermonde determinant~${\sqrt{\Delta(F)}}$
as a common denominator
in the sense that~${\sqrt{\Delta(F)}}c_i$
is an algebraic integer for each~$i=1,\ldots,k$. Note that $Mc_i$ is a multiple of the algebraic integer ${\sqrt{\Delta({\mathbb K})}}$, which is an algebraic integer in ${\mathcal O}_{\mathbb K}$.

\subsection{The algebraic Dold condition}

We say that a sequence of algebraic
integers~$(v_n)_{n\ge 1}$ satisfies the
\emph{algebraic Dold condition} if~(D)
is satisfied as algebraic integers. That is, if
$$
\frac{1}{n}\sum_{d\vert n} \mu(n/d)v_d\in {\mathcal O}_{\mathbb K}
$$
for all~$n\ge1$.
Our strategy is to find~$s$
such that~$(\lambda_i^{ns})_{n\ge 1}$ satisfies the algebraic
Dold condition for~$i=1,\ldots,k$.
Since linear combinations with algebraic integer
coefficients of sequences which satisfy
the algebraic Dold condition also satisfy the
algebraic Dold condition, this allows us
to deduce that
$$
M\frac{1}{n} \sum_{d\vert n} \mu(n/d)u_{d^s}=M\sum_{i=1}^k c_i\frac{1}{n}\sum_{d\vert n} \mu(n/d)\lambda_i^{d^s}
$$
is both a rational number and an algebraic integer, so an integer,
verifying~(D).

Fix~$\lambda:=\lambda_i$ for some~$i=1,\ldots,k$.
We need to find out when
\begin{equation}\label{eq:algDoldargument1}
\frac{\sqrt{\Delta({\mathbb K})}}{n} \sum_{d\vert n} \mu(n/d) \lambda^{d^s}
=
\frac{{\sqrt{\Delta({\mathbb K})}}}{n} \sum_{d\vert n} \mu(d) \lambda^{(n/d)^s}
\end{equation}
is an algebraic integer for all~$n\ge1$.
We write~$m:=\prod_{p\vert n} p={\text{\rm rad}}(n)$ for the radical of~$n$.
Changing the order of summation to
complementary divisors as shown and
restricting to squarefree numbers
in the summation on the
right-hand side (as the M\"obius function
vanishes on all other terms), the numerator
in~\eqref{eq:algDoldargument1} is
\begin{equation}\label{eq:algDoldargument2}
S:={\sqrt{\Delta({\mathbb K})}}\sum_{d\vert m} \mu(d) \lambda^{(n/d)^s}.
\end{equation}
Clearly~$S$ is a multiple of~$1$,
so we may assume that~$n$ (and hence~$m$) exceeds~$1$.
Let~$p$ be a prime divisor of~$n$ and let~$w$ be the exact exponent of~$p$ in~$n$,
written~$p^w\| n$.
Writing as usual~$\omega(n)$ for the number
of distinct prime divisors of~$n$, the sum~\eqref{eq:algDoldargument2}
therefore
has~$2^{\omega(n)}$ divisors, half of which are multiples of $p$ and half of which are not. Thus,
the above sum can be grouped into~$2^{\omega(n)-1}$
pairs indexed~$(d,p)$, where~$d$ is a divisor of~$\frac{m}{p}$,
giving
\begin{align*}
S & =  {\sqrt{\Delta({\mathbb K})}}\sum_{d\vert\frac{m}{p}} \left(\mu(d) \lambda^{(n/d)^s}+\mu(pd) \lambda^{n/(dp)}\right)\\
& =  {\sqrt{\Delta({\mathbb K})}}\sum_{d\vert\frac{m}{p}}\pm \lambda^{(n/(pd))^s}\left(\lambda^{(n/pd)^s(p^s-1)}-1\right).
\end{align*}
Thus,
it is sufficient to show that if~$p^w\| n$, then~$\frac{S}{p^w}$ is an algebraic integer. We let $\pi$ be a prime ideal divisor of $N$
in ${\mathbb K}$ with $\pi^e\| p$ and $N_{{\mathbb K}/{\mathbb Q}}(\pi)=p^f$. We put
$$
A_d:=\lambda^{(n/dp)^s}\qquad {\text{\rm and}}\qquad B_d:=\lambda^{(n/pd)^s(p^s-1)}-1,
$$
and observe that our aim is to find a condition for $s$ divisible by $e(G)$ such that one would have
$$
\nu_{\pi}({\sqrt{\Delta({\mathbb K})}}A_dB_d)\ge ew.
$$
Here, $\nu_{\pi}(\alpha)$ is the exponent of the prime ideal $\pi$ in the factorization of $\alpha{\mathcal O}_{\mathbb K}$.
Observe that the Different Theorem implies that
$$
\nu_{\pi}(\Delta({\mathbb K}))\ge fe(e-1)\ge e(e-1).
$$

\subsubsection{$\pi\vert\lambda$}
\label{sec:1}

In this case~$\nu_{\pi}(B_d)=0$ and
$$
\nu_{\pi}(A_d)\ge \left(\frac{n}{pd}\right)^s\ge p^{s(w-1)}\ge 2^{s(w-1)}.
$$
We need to check that
$$
2^{s(w-1)}+e(e-1)/2\ge ew.
$$
This is clear when~$w=1$.
since then the inequality to be proved becomes~$1+e(e-1)/2\ge e$ which holds for all~$e\ge 1$.
This is also clear when~$e=1$ since then it is implied by~$2^{s(w-1)}\ge 2^{w-1}\ge w$.
Finally, if~$e\ge 2,~w\ge 2$, then~$s\ge N\ge e\ge 2$, so~$sw\ge 4$. Since~$w-1\ge w/2$, it suffices to show that~$2^{sw/2}\ge sw$, which is equivalent to~$2^{sw}\ge (sw)^2$,
which holds since~$sw\ge 4$.

\subsubsection{$p\!\nmid\!\lambda$}
\label{sec:2}

In this case,~$\nu_{\pi}(A_d)\ne 0$.
Write~$(n/pd)^s=\alpha p^{s(w-1)}$, and~$p^s-1=\beta(p^f-1)$,
where~$\pi\!\nmid\!\alpha\beta$. This last formula holds since~$f\vert e(G)\vert s$.
The analogue of Euler's theorem for number fields implies that
$$
\lambda^{p^{s(w-1)}(p^f-1)}\equiv 1\pmod {\pi^{s(w-1)+1}}.
$$
Thus
$$
\nu_{\pi}(B_d)\ge s(w-1)+1.
$$
So it suffices to verify that
$$
s(w-1)+1+e(e-1)/2\ge ew.
$$
This is clear if~$w=1$ since then the left-hand side is~$1+e(e-1)/2\ge e$. It is also clear if~$e=1$, since then the left-hand side is~$s(w-1)+1\ge (w-1)+1=w$. Thus, we assume that~$e\ge 2$ and~$w\ge 2$. Since~${\mathbb K}$ is Galois, we have that~$e\vert N$ and~$s\ge N$. If~$s\ge 2e$, then it suffices to show that
$$
2e(w-1)+1\ge ew.
$$
This is equivalent to~$ew-2e+1\ge 0$, which holds since~$w\ge 2$. Finally, if~$2e>s\ge N$,
we have~$e>N/2$ and~$e$ is a divisor of~$N$ so~$s=e=N$. So
we need to show that
$$
N(w-1)+1+N(N-1)/2\ge Nw,
$$
which is equivalent to~$1+N(N-1)2\ge N$, which is clear for any~$N\ge 2$.

\subsection{The sign condition}
\label{sectionSignCondition}

We still need to deal with the sign condition~(S), for
which we may use the following observation from~\cite{Pu}:
It is sufficient
to show that~$u_{(2n)^s}\ge n u_{n^s}$
for all~$n\ge 1$.
To see this, let~$\lambda$ be a real root
larger than~$1$ of~$F(X)=0$.
This exists by the intermediate
value theorem, since the hypotheses on
the coefficients~$a_1,\dots,a_k$ show that~$F(1)<0$, and~$F(x)\to\infty$
as~$x\to\infty$.
Note that~$u_n\ge \lambda^{n-k}$ always holds,
again by the hypotheses on the coefficients.
Indeed, it holds for~$n=1,\ldots,j$
because in this range~$u_j\ge 1\ge \lambda^{j-k}$,
and so it holds for all~$n\ge1$ by induction
since~$a_i\ge 0$ for~$i=1,\ldots,k$.
Moreover,~$u_n\le \lambda^{n+n_0}$
for~$n_0\ge\lceil \log\max\{u_1,\ldots,u_k\}/\log \lambda\rceil$.
Again this inequality holds for~$n=1,\ldots,k$,
so it will hold for all~$n\ge 1$
by induction.
Armed with these estimates, we now need to show that
$$
\lambda^{(2n)^s-k}\ge n\lambda^{n^s+n_0}
$$
or, equivalently,
\begin{equation}
\label{eq:n0}
n^s(2^s-1)\ge k+n_0+\frac{\log n}{\log\lambda}.
\end{equation}
To see this, first let~$n_1\ge 2(n_0+k)$
satisfy~$\frac{\log n}{\log\lambda}\le\frac{n}{2}$
for~$n\ge n_1$.
Then for~$n\ge n_1$
we have~$k+n_0\le\frac{n_1}{2}\le\frac{n}{2}$,
so the right-hand side of~\eqref{eq:n0} is
at most~$n$.
It follows that~\eqref{eq:n0} holds, since~$n(2^s-1)\ge n$
is clear.
For~$n\le n_1$
the right-hand side is at most~$n_0+k+\frac{\log n_1}{\log\lambda}$
and the left-hand side is at least~$2^s(2^s-1)\ge 2^{e(G)\ell}(2^{e(G)\ell}-1)$ and this is larger
than~$n_0+k+\frac{\log n_1}{\log\lambda}$
once~$\ell\ge\ell_0$. Thus,
if~$s$ is a sufficiently large multiple of~$e(G)$,
then the sign condition~(S) holds.

\section{The proof of Theorem \ref{thm:2}}

For the particular case of the~$k$-generalized Fibonacci sequence,
it is well-known that the associated
characteristic polynomial
$$
F^{(k)}(X)=X^k-X^{k-1}-\cdots-1
$$
has simple zeros (see Miles~\cite{MR123521}, for example), and
that the largest real zero~$\lambda^{(k)}$
is increasing in~$k\ge2$ and has~$\lambda^{(k)}\to2$
as~$k\to\infty$.
In particular,
writing
$$
\lambda=\lambda^{(k)}\ge\lambda^{(2)}=\frac{1+\sqrt{5}}{2}
$$
we have~$F_k^{(k)}<2^k<\lambda^{2k}$, so we can take $n_0=2k$
in the notation of Section~\ref{sectionSignCondition}.
Thus,~$n_1\ge 2(n_0+k)=6k$
must be such that~$\frac{\log n}{\log\lambda}\le\frac{n}{2}$,
which is implied by~$6\log n\le n$, which
certainly
holds for~$n\ge
10k\ge 20$. So we can take~$n_1=10k$.
Consequently, for~$n\le n_1$, the right-hand
side in~\eqref{eq:n0} is at most~$10k$ and the left-hand side
is at least~$n^s(2^s-1)\ge n^{N_k}(2^{N_k}-1)\ge 10k$
for all~$n\ge1$ and~$k\ge~3$
where~$N_k=k!$. For~$k=2$,
the above inequality
fails for~$n=1,2$, but in these
cases~$F_{(2n)^2}\ge F_4=3> 2F_2$ holds anyway.
Hence, we can take~$\ell_0=1$ for any~$k\ge 2$.

\section{The proof of Theorem \ref{thm:3}}

\subsection{The case~$Q=0$}
In this
degenerate case we may take two approaches
(for~$P>0$ at any rate), and we include both to
illustrate the two points of view.

\noindent{{\bf{Arithmetic proof:}}
We have~$u_n=P^{n-1}$, so~$c_1=1/P$
in Binet's formula and~$\lambda=P$.
Going through the proof of Theorem~\ref{thm:1},
we see that we need $P\vert M$ to deal with the denominator of~$c_1$.
Next, in case~$p$ does not divide~$P$, we are in the case
from Section~\ref{sec:2}, and then~$p^w\vert S$
whenever~$p^w\vert n$. In the case of Section~\ref{sec:1}, we have that
if~$p\vert P$, then~$e=N=1$. We saw in the proof of Theorem~\ref{thm:1}
for this case that if~$p\vert P$ and~$w\ge 2$, then $\nu_p(A_d)\ge 2^{w-1}>w$, whereas for~$w=1$,
we have $\nu_p(A_d)]\ge 1=w.$ So, in fact we even
see that~$(|P|u_{n^2})_{n\ge 1}$
satisfies the Dold condition in this case, and we do not need the factor~$P^2$.

\noindent{\bf{Dynamical proof for $P>0$:}}
Here~$(Pu_n)_{n\ge1}$
is the sequence~$(P,P^2,\dots)$,
which we identify as~$({\text{\rm Fix}}(T^n))_{n\ge1}$,
where~$T\colon X\to X$ is the shift map
on the \emph{full~$P$-shift}~$X=\{1,2,\dots,P\}^{\mathbb{Z}}$.
Taking the union of~$P$ disjoint copies of
this system produces a map~$S$ with~$({\text{\rm Fix}}(S^n))_{n\ge1}
=(P^2,P^3,P^4,\dots)$. It follows that~$(P^2u_n)_{n\ge1}$
is realizable, and in particular satisfies~(D).
On the other hand, sampling along a monomial subsequence
always preserves
realizability (no other polynomials have this property
by~\cite{MR4002553})
so~$(P^2u_{n^2})$ is also realizable.

\subsection{The case $(P,Q)=(0,1)$}
This is an excluded case of
Theorem~\ref{thm:1}. In this case,~$u_n=0$
if~$n$ is even and~$u_n=1$ if~$n$ is odd.
So, in~(D), the sum~$S$ is zero if~$n$ is odd.
Further, in the sum~$S$, for every prime~$p$ dividing~$m$,
the amounts~$(n/d)^2$ and~$(n/(pd))^2$ are both even or both
odd, so this difference is zero
unless~$p=2$ and one of~$(n/d)^2$ and~$(n/(2d))^2$
is even and the other is odd.
But the only chance for this to happen is when~$2\| n$,
and in this last case the prime~$2$ from the denominator of the Dold ratio can be absorbed into~$|\Delta(F)|=|P^2-4Q|=4$.

\subsection{Remaining cases}
The cases~$(P,Q)\ne (P,0),~(0,1)$ follow from Theorem \ref{thm:1}.



\begin{thebibliography}{10}

\bibitem{MR1990179}
G.~Everest, A.~van~der Poorten, I.~Shparlinski, and T.~Ward, \emph{Recurrence
  sequences}, in \emph{Mathematical Surveys and Monographs} \textbf{104}
  (American Mathematical Society, Providence, RI, 2003).
\newblock \verb|https://doi.org/10.1090/surv/104|.

\bibitem{MR1555035}
P.~Fatou,  `S\'eries trigonom\'etriques et s\'eries de {T}aylor', \emph{Acta
  Math.} \textbf{30} (1906), no.~1, 335--400.
  \newblock\verb|https://doi.org/10.1007/BF02418579|.

\bibitem{MR4002553}
S.~Jaidee, P.~Moss, and T.~Ward,  `Time-changes preserving zeta functions',
  \emph{Proc. Amer. Math. Soc.} \textbf{147} (2019), no.~10, 4425--4438.
\newblock \verb|https://doi.org/10.1090/proc/14574|.

\bibitem{MR2043329}
P.~A. Martin,  `The {G}alois group of {$x^n-x^{n-1}-\dots-x-1$}', \emph{J. Pure
  Appl. Algebra} \textbf{190} (2004), no.~1-3, 213--223.
\newblock \verb|https://doi.org/10.1016/j.jpaa.2003.10.028|.

\bibitem{MR123521}
E.~P. Miles, Jr.,  `Generalized {F}ibonacci numbers and associated matrices',
  \emph{Amer. Math. Monthly} \textbf{67} (1960), 745--752.
\newblock \verb|https://doi.org/10.2307/2308649|.

\bibitem{MR3195758}
G.~T. Minton,  `Linear recurrence sequences satisfying congruence conditions',
  \emph{Proc. Amer. Math. Soc.} \textbf{142} (2014), no.~7, 2337--2352.
\newblock \verb|https://doi.org/10.1090/S0002-9939-2014-12168-X|.

\bibitem{pm}
P.~Moss, \emph{The arithmetic of realizable sequences} (Ph.D. thesis,
  University of East Anglia, 2003).

\bibitem{MT}
P.~Moss and T.~Ward,  `Fibonacci along even powers is (almost) realizable',
  \emph{Fibonacci Quart.} \textbf{60} (2022), no.~1, 40--47.
\newblock \verb|https://www.fq.math.ca/Abstracts/60-1/moss.pdf|.

\bibitem{MR2090972}
M.~R. Murty and J.~Esmonde, \emph{Problems in algebraic number theory}, in
  \emph{Graduate Texts in Mathematics} \textbf{190} (Springer-Verlag, New York,
  second ed., 2005).

\bibitem{Pu}
Y.~Puri, \emph{Arithmetic of numbers of periodic points} (Ph.D. thesis,
  University of East Anglia, 2001).

\bibitem{MR1873399}
Y.~Puri and T.~Ward,  `Arithmetic and growth of periodic orbits', \emph{J.
  Integer Seq.} \textbf{4} (2001), no.~2, Article 01.2.1, 18.
\newblock \verb|https://www.emis.de/journals/JIS/VOL4/WARD/short.pdf|.

\bibitem{MR1866354}
Y.~Puri and T.~Ward,  `A dynamical property unique to the {L}ucas sequence',
  \emph{Fibonacci Quart.} \textbf{39} (2001), no.~5, 398--402.
\newblock \verb|https://www.mathstat.dal.ca/FQ/Scanned/39-5/puri.pdf|.

\bibitem{MR1622060}
D.~A. Wolfram,  `Solving generalized {F}ibonacci recurrences', \emph{Fibonacci
  Quart.} \textbf{36} (1998), no.~2, 129--145.
\newblock \verb|https://www.fq.math.ca/Scanned/36-2/wolfram.pdf|.

\end{thebibliography}

\providecommand{\bysame}{\leavevmode\hbox to3em{\hrulefill}\thinspace}

\end{document}